\documentclass[final]{siamltex}
\usepackage{graphicx,amssymb}

\def\LU{L\kern .8pt U}\def\PA{P\kern -.9pt A}
\def\PtA{\tilde P\kern -.9pt A}
\def\qed{~~\hbox{\vrule width 3.5pt depth 5 pt height 2.5 pt}}

\def\pput(#1,#2)#3{\noindent\smash{\raise#2pt\hbox to 0pt
   {\kern #1pt #3\hss}}\ignorespaces}

\title{Instability of Gaussian elimination is exponentially rare (proof of partial result)}

\author{Lloyd N. Trefethen\thanks{School of Engineering and Applied Sciences,
Harvard University, Cambridge, MA 02139, USA ({\tt trefethen@seas.harvard.edu}).}}

\begin{document}

\maketitle

\begin{abstract}
Gaussian elimination with partial pivoting is the standard
algorithm for solving a dense $n\times n$ system of linear
equations $Ax=b$.  The ``growth factor'' $\rho\,$ for this
process may be as large as $2^{n-1}$, and when $\rho$ is
large, the algorithm is unstable.  Nevertheless, long
experience has shown that Gaussian elimination is resoundingly
stable in practice.  Experiments
indicate that among random matrices with independent normally
distributed entries, growth factors $\rho \gg n^{1/2}$ are
in a precise sense exponentially rare.
Here we take a step toward confirming this observation by proving
it not for the full growth factor, which depends
on all the entries of the upper-triangular factor $U$, but for the 
corner entry $u_{nn}^{}$.
\end{abstract}

\begin{keywords}
Gaussian elimination, growth factor, stability
\end{keywords}

\pagestyle{myheadings}
\thispagestyle{plain}

\markboth{L. N. TREFETHEN}
{STABILITY OF GAUSSIAN ELIMINATION}

\section*{Preliminary note}  This paper was essentially written in
2015; just a few updates have been added in October 2026.  I have decided
to make this text available at the same time as John Urschel of MIT posts
a manuscript in which he apparently settles our conjecture of
section 5 in the affirmative,
using quite different arguments~\cite{urschel}.

\section{Introduction}

For six decades, Gaussian elimination with partial pivoting
(henceforth GE) has been the standard algorithm for solving
general systems of linear equations $Ax=b$.  Trillions of systems
have been solved by GE over the years, and dimensions $n$ in the
tens of thousands
are now routine.  In all these computations,
numerical instabilities have virtually never been
reported.\footnote{The three reports I am aware of are by
Foster~\cite{fos}, Liu and Russell~\cite{rus}, and Wright~\cite{wri},
but these cases seem unusual enough not to have had any
impact on the reliance on GE in practice.}

This almost perfect reliability of GE has not been explained.  In the
1940\kern .4pt s, it was thought that algorithms
of this type would be unstable.  Hotelling predicted exponential
amplification of rounding errors in 1943~\cite{hot}, as did
Bargmann, Montgomery, and von Neumann in 1946~\cite{bmv}. Gaussian
elimination was then investigated by Wilkinson, who in a celebrated
analysis in 1961 showed that the problem of stability could be
reduced to the problem of bounding the {\it growth factor\/}
$\rho\kern .7pt$ of the elimination process~\cite{wil61}: the ratio of the
maximum element generated during the elimination to the maximum
element of the original matrix.  Wilkinson proved that GE can
be unstable only if $\rho$ is large.  This classical
analysis appears in numerous textbooks, including the standard
reference on matrix computations by Golub and Van Loan~\cite{gvl}.

This would be the end of the matter if $\rho$ were always small.
In fact, $\rho$ can be as large as $2^{n-1}$.  Examples have
been known from the beginning that achieve this bound,
and for such examples, instability of magnitude $O(2^n)$ is
indeed observed, making GE useless for almost all values
of $n$.\footnote{For example, let $A$ be the $50\times 50$
matrix with $1$ on the diagonal, 1 in the rightmost column, and
$-0.9$ in each position below the diagonal.  About 12 digits
of accuracy are lost to rounding errors when $Ax=b$
is solved by GE in standard 16-digit arithmetic for a typical
right-hand side $b$, whereas a solution by QR factorization or
singular value decomposition gives close to machine precision.
The disparity persists if the entries of $A$ are modified
by perturbations of size $10^{-14}$.}
Thus the behavior of GE is two-headed: hopelessly
unstable for certain problems, abundantly stable in practice.

\begin{figure}[h]
\begin{center}
\includegraphics[scale=.7]{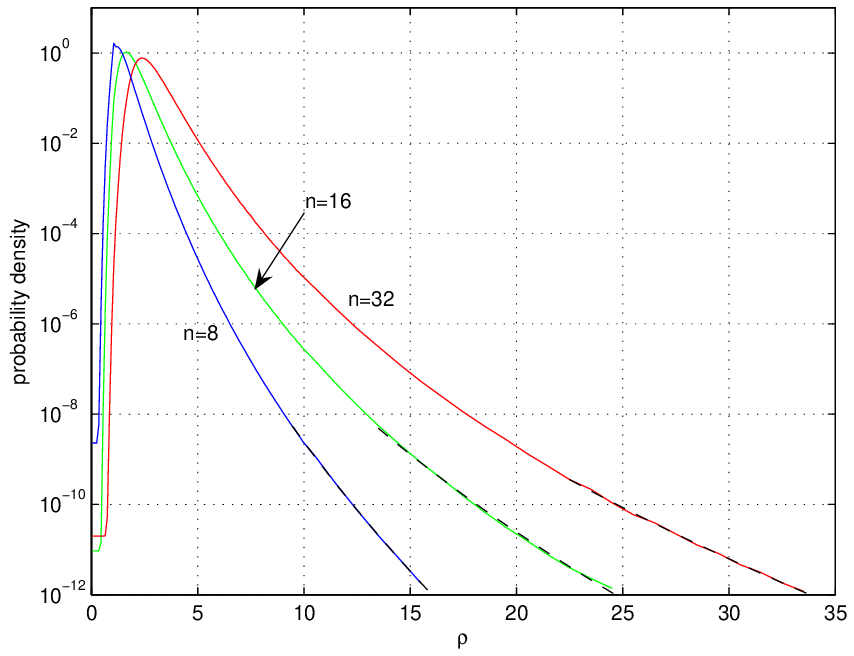}
\end{center}
\caption{\label{fig-DM}Figure\/ $8$ of\/ $\hbox{\rm \cite{dm}}$,
showing empirically that
large growth factors $\rho$ are extremely rare among random
matrices of dimensions $n = 8$, $16$, and\/ $32$;
see also $\hbox{\rm \cite{tb}}$ and $\hbox{\rm \cite{ts}}$.
(The definition of $\rho$ here is slightly different
from that of\/ $(\ref{rho})$ in that the numerator
is restricted to $k=n$.)}
\end{figure}

In view of the central position of GE in scientific computing,
one might expect that much effort would have been devoted to
explaining this anomaly.  However, it is not so.  Perhaps because
Wilkinson's introduction of backward error analysis for this
and other problems of numerical stability was itself such
a striking advance,
or perhaps because problems of a probabilistic nature are not
always appealing to numerical analysts, the investigation of the
stability of GE has remained essentially where Wilkinson left
it in 1961.  The main discussion of the problem that has been
published to date is Chapter 22 of \cite{tb}, which the reader of
this article is encouraged to consult for background along with
the earlier related paper~\cite{ts}.\ \ Figure~\ref{fig-DM}
reproduces a figure by Driscoll and Maki 
from \cite{dm}, an improvement of Figure~22.8
of~\cite{tb}, showing the extreme rarity of even modestly large
values of $\rho$ among random matrices.\footnote{In the words of
Driscoll and Maki~\cite{dm}: ``Our results suggest, for example,
that the occurrence of an $8\times 8$ matrix with a growth factor
of $40$ is on the order of a once-in-the-age-of-the-universe
event.''} However, no theorems have been published.\footnote{The
closest relevant results are those of the
{\em smoothed analysis\/} of Sankar, Spielman, and
Teng, which considers the behavior of
Gaussian elimination in a cloud of matrices that
are perturbations of an underlying matrix.  However, the main focus
of their work is on behavior of elimination without
pivoting~\cite{sst}.  A result on elimination with
pivoting can be found in Sankar's unpublished thesis~\cite{san}, showing
that growth factors only rarely exceed a
quantity on the order of $n^{\log n}$.}

Here we take a step toward proving the theorem that the experiments have
long suggested must be true (and for which a $\$1000$ prize
was offered in~\cite{tre}).  We show that among random
matrices with independent standard normally distributed entries,
corner entries $|u_{nn}^{}| \gg n^{1/2}$ are exponentially rare in
a precise sense that we define.  The as-yet-unproved theorem, listed
as the Conjecture in section 5, is the
assertion that the same is true for all the entries $u_{ij}^{}$.
The crux of our argument is to prove that whereas random
triangular matrices are exponentially ill-conditioned~\cite{vt}, the
triangular factors of random matrices computed by GE are
exponentially unlikely to be ill-conditioned at all.\footnote{See
Figure~22.3 of~\cite{tb}.}

\section{Problem setting and the definition of ``exponentially unlikely''}

In this article we employ the following notions of
random vectors and matrices (called the {\em Gaussian
ensemble\/} in the random matrix literature).

\medskip

{\em Definition.}  A {\em random $n$-vector\/} is an $n$-vector
whose entries are independent samples from the standard normal
distribution $N(0,1)$.  A {\em random\ $n\times n$ matrix\/}
is an $n\times n$ matrix whose entries are independent samples
from this distribution.

\medskip

It is well known that the probability distribution of the
$k$-dimensional spaces spanned by $k$ random $n$-vectors (with
probability 1) is isotropic in ${\mathbb R}^n$, invariant under orthogonal
transformations; likewise the $(n-k)$-dimensional space
orthogonal to $k$ random $n$-vectors is isotropically oriented.

Gaussian elimination with partial pivoting consists of a
sequence of elimination steps $k = 1,2,\dots, n-1$ in which
at step $k$, row $k$ is first interchanged (``pivoted'')
with a later row if necessary so that the largest entry
in or below the diagonal in column $k$ (``the pivot'') is
moved to the diagonal,\footnote{Here and throughout, when we
apply terms like ``large'' or ``small'' to a number, we are
speaking of its absolute value.} and then multiples of row
$k$ are subtracted from rows $k+1,\dots,n$ so that column $k$
becomes zero below the diagonal.\footnote{The term ``partial''
refers to searching for a largest entry in a single column.
A search over all columns is {\em complete pivoting,}
but this is rarely used since partial pivoting is so effective
in practice.  If no pivoting at all is employed, Gaussian
elimination is unstable for general matrices when it does
not break down entirely: both inaccurate in
practice and provably unstable in theory, with unbounded growth
factors whose expected values in a probabilistic analysis are
infinite.} Ties in the choice of pivot are broken arbitrarily;
the treatment of ties will be of no consequence for us since
they occur with probability zero.  The result at the end of
elimination is a factorization $\PA = \LU$
in which $L$ is lower-triangular with $1$ on the diagonal and
$|\ell_{ij}|\le 1$ below the diagonal, $U$ is upper-triangular,
and $P$ is a permutation matrix, meaning that each row and each
column contains a single nonzero entry~$1$.\ \ Full details on
GE can be found in \cite{gvl}, \cite{ste}, or \cite{tb}.

The {\em growth factor\/} for a particular
matrix $A$ is defined by
\begin{equation}
\rho = {\max_{i,j,k} |u_{ij}^{(k)}| \over \max_{i,j} |a_{ij}|}. 
\label{rho}
\end{equation}
Here $u_{ij}^{(k)}$ denotes the $(i,j)$ entry of $U^{(k)}$, the
partly triangular matrix that has been constructed
after $k$ steps of elimination, zero below the diagonal in the
first $k$ columns.

The following definition is employed throughout this paper
to enable us to present the argument in a
conversational style while still being mathematically precise.

\medskip

{\em Definition.} Let $f$ and $g$ be quantities associated
with random $n\times n$ matrices
for all sufficiently large $n$.  We say that it is\/
{\em exponentially unlikely that $f\gg g$},
or in words, {\em exponentially unlikely that $f$ is
much greater than $g$}, if for any $\alpha >0$ and any $M>0$,
\begin{equation}
\hbox{\rm probability} \left( |f(n)| > |g(n)| n^\alpha \right) < n^{-M}
\label{highly}
\end{equation}
for all sufficiently large $n$.

\section{First theorem: \boldmath$n^{3/2}$ for the final entry of\ $U$}

Our main theorem, in the next section, establishes an
$O(n^{1/2})$ bound on the corner entry $u_{nn}^{}$ of the matrix $U$.  In this section,
we give a first part of the argument that is enough to prove a
weaker $O(n^{3/2})$ bound (or $O(n)$ for complex matrices).

\medskip

{\em {\bf Theorem 1.}  Let GE be applied to a random
$n\times n$ matrix.  It is exponentially unlikely that
$u_{nn}^{}$ is much greater than $n^{3/2}$.}

\medskip

{\em Proof.} GE produces a factorization $\PA = \LU$, where
each subdiagonal entry of $L$ is no greater than $1$.
By induction, it can be shown that this implies that the entries in
columns $n, n-1, n-2, \dots$ of $L^{-1}$ are no
greater than $1, 1, 2, 4, 8, \dots.$

\begin{figure}
\begin{center}
\vskip .3cm
\includegraphics[scale=.65]{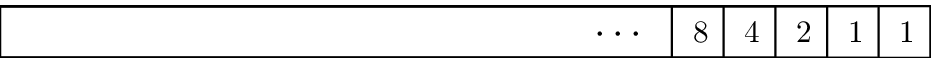}
\end{center}
\caption{\label{fig-thm1}Sketch of the key step of the proof of
Theorem\/~$1$.  The last row of $L^{-1}$, $m$, is a random\/
$n$-vector\/ $w$, ``diabolically pivoted,'' then normalized to have
last entry $1$, and its final entries are bounded
by $1,1,2,4,8,\dots.$  If the largest entry of\/ $m$ is of
size $n\kern .4pt \varepsilon^{-1}$ with $\varepsilon \ll 1$,
such a configuration
requires the vector $w$ to have $O(\log \varepsilon)$ entries
of size $\ll \varepsilon n^{-1}$, which is
exponentially unlikely.}
\end{figure}

Let $m$ denote the last row of $L^{-1}$.  The equation $\PA=\LU$
implies $U = L^{-1} \PA$.  Thus $u_{nn}^{}$ is the product
of $m$ with the last column of $\PA$, which consists of $n$
independent samples from $N(0,1)$, unaffected by pivoting since
the elimination process has not made use of data from this
column.
This implies that $u_{nn}^{}$ is a normally distributed
random variable of mean zero whose standard deviation is no greater
than $n^{1/2} \max_j |m_j|$, and therefore it
is exponentially unlikely that $u_{nn}^{}$ exceeds 
$\max_j |m_j|$ by much more than a factor of $n^{1/2}$.  Thus to
prove Theorem 1, it is enough to show that it is exponentially
unlikely that any entry of $m$ is much greater than~$n$.

In the remainder of the argument, we make no careful analysis
of pivoting.  On the contrary, we assume that matrices may be
``diabolically pivoted,'' that is, we derive probabilistic
bounds that apply even in the worst case in which the pivoting
happens to be as adverse as possible.  Nevertheless, pivoting
enters the analysis in a crucial way in implying that the final
entries of $m$ are no greater than $1,1,2,4,8,\dots.$

Now since the last row of $U$ is zero except in the $(n,n)$
position, $m$ is orthogonal to the first $n-1$ columns of $\PA$.
These columns span a random $(n-1)$-dimensional space, apart
from the pivoting.  Therefore $m$, being orthogonal to these,
is a permuted and then normalized random $n$-vector.\footnote{For
Gaussian elimination without pivoting, $m$ is not permuted---it
is simply $m = w/w(n)$, where $w$ is a random vector.  This
observation immediately makes probabilistic analysis of the
entries of $L^{-1}$ and consequently $U$ straightforward for
elimination without pivoting.} We can describe its structure
like this:

\medskip 

\begin{enumerate}
\item
Let $w$ be a random $n$-vector, with entries arbitrarily permuted;
\item 
$m = w/w(n)$.
\end{enumerate}

\medskip 

If $w$ is a random $n$-vector, its smallest entries scale as
$n^{-1}$ and its largest as $1$.  The probability that
$w$ has an entry of size $\varepsilon\kern .4pt n^{-1}$ is $O(\varepsilon)$
as $\varepsilon \to 0\kern .7pt $; the probability that
it has two such entries
is $O(\varepsilon^2)$, and so on.  (For complex
vectors, these probabilities apply
to entries of size $\varepsilon^{1/2}n^{-1/2}$.)
For its pivoted normalized
vector $m$ to have an entry greater
than $n\kern .4pt \varepsilon^{-1}$ with
$\varepsilon \ll 1$, $w$ would
need to have $O(\log \varepsilon)$ such small entries because of
the $1,1,2,4,8,\dots$ constraint.  To be specific, suppose
we define
\begin{equation}
\mu = \hbox{\rm round}({\textstyle{1\over 2}}|\log_2 \varepsilon\kern .5pt |).
\label{mudef}
\end{equation}
Then $w$ would have to have at least $\mu$ entries smaller
than $\varepsilon^{1/2} n^{-1}$, and for $\varepsilon\ll 1$, this
is exponentially unlikely.\footnote{Discussion
with a reader suggests the following analogy.
Suppose somebody has tossed 100 coins
and ordered the results in a diabolical fashion you know
nothing about.  If the first 50 results are all heads, this cannot
count as surprising, because it could be due to
the ordering.  If the first 75 results are all heads, however,
that is surprising indeed, for even diabolical ordering would be
unlikely to be able to achieve that with 100 coin tosses.}\qed

\section{Second theorem: \boldmath$n^{1/2}$ for the final entry of\ $U$}

Our second theorem improves the bound just
established from $n^{3/2}$ to $n^{1/2}$.  

\medskip

{\em {\bf Theorem 2.}  Let GE be applied to a random
$n\times n$ matrix.  It is exponentially unlikely that
$u_{nn}^{}$ is much greater than $n^{1/2}$.}

\medskip

{\em Proof.}  
We repeat the argument of Theorem 1, but now, we improve the
bound on the entries of $m$, the final row of $L^{-1}$, from $n$
to $1$.  Specifically, to establish Theorem~2, it is enough to
prove---as seems clear in numerical experiments---that it is
exponentially unlikely that any entry of $m$ is much greater
than~$1$.

Suppose $m$ contains an entry of size $\varepsilon^{-1}$ for
some $\varepsilon \ll 1$.  The key idea is to consider the
configuration $\mu$ steps before the end of elimination, with
$\mu$ defined by (\ref{mudef}).  At this stage, the elimination
has not yet made use of data from the last $\mu$ columns of $A$.
Therefore, we can apply a random orthogonal operation to these
last $\mu$ columns without changing the probability distribution.
From this point on, data from these columns will begin to be
brought in, and further row interchanges will take place, but
all these row interchanges will be confined to the final $\mu$
rows of~$A$.

\begin{figure}
\begin{center}
\vskip .1cm
\includegraphics[scale=.65]{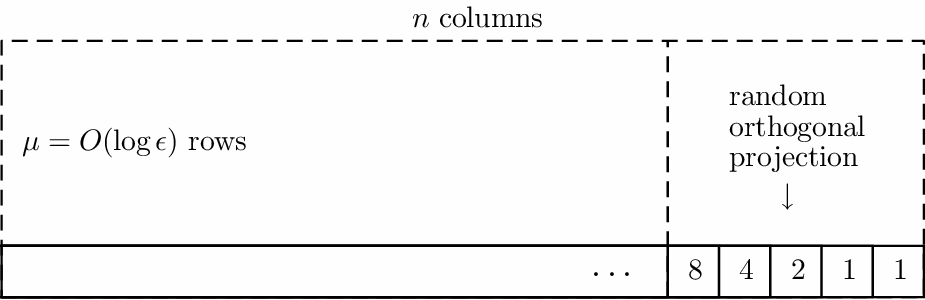}
\end{center}
\caption{\label{fig-thm2}Sharpening of the
key step of the proof of Theorem\/~$1$ to establish
Theorem\/~$2$.  At a stage $\mu = O(\log \varepsilon)$
steps from the end of elimination, we know that all further pivoting
will affect only the final $\mu$ columns of $L^{-1}$.  Therefore,
for the final vector $m$ to have a reasonable chance of
having an entry much larger
than $1$, an entire $\mu\times \mu$
block of a diabolically pivoted random $\mu\times n$ matrix would have
to be small or at least close to low
rank.  This is exponentially unlikely.}
\end{figure}

This enables us to set down a sharper characterization of $m$
than before, sketched in Figure~\ref{fig-thm2}.  Note that like
the vector $w$ of the last section, the matrix $W$ defined in
line 1 below is not one that appears explicitly during GE.\ \
Rather, it is an artifice that spans subspaces of
${\mathbb R}^n$ with the appropriate probability distribution.
Line 2 corresponds to applying a random orthogonal operation
to the $\mu$-dimensional space still in play, before selecting
a particular representative vector.

\medskip

\begin{enumerate}
\item
Let $W$ be a random $\mu\times n$ matrix, with columns arbitrarily permuted;
\item 
Let $w_0 = q^TW$, where $q$ is a random $\mu$-vector normalized
to $\|q\|=1$;
\item
Let $w$ be a further arbitrary permutation of $w_0$, in its
final $\mu$ entries only;
\item
$m = w/w(n)$.
\end{enumerate}

\medskip

\noindent
As before, $m$ has the additional property that its final
entries are bounded by $1,1,2,$ $4,8,\dots .$ Specifically,
if $\mu$ is defined by (\ref{mudef}), the final $\mu$
entries of $m$ must all be bounded by $\varepsilon^{-1/2}$.
Now this condition is unaffected by step 3 in the list above,
a permutation of just the final $\mu$ entries of $w_0$.  We are
thus done if we can prove the following assertion extracted from
steps 1, 2, and 4: if $W$ is a random $\mu\times n$ matrix,
it is exponentially unlikely that $W$ has a $\mu\times \mu$
submatrix $W_\mu$ (i.e., a selection of some $\mu$ distinct
columns of $W$) with the property that $q^TW_\mu$, where $q$
is a random $\mu$-vector normalized by $\|q\|=1$, has all its
entries smaller than $\varepsilon^{1/2}$ for $\varepsilon \ll 1$.

Intuitively speaking, here is why such a configuration is
exponentially unlikely.  Though $W$ has some small entries (of
order $n^{-1}$), it is exponentially unlikely to have any small
{\em columns,} i.e., columns all of whose entries are small.
So any matrix $W_\mu$ will consist of columns that are not small.
It is algebraically unlikely that the inner product of a fixed
random $q$ with such a column is small, and one would accordingly
expect it to be exponentially unlikely that its inner product
with all $\mu$ of the columns is small.  This last expectation,
however, cannot immediately be justified, since the columns of
$W_\mu$ are not independent due to the pivoting.

To make a rigorous argument, we proceed as follows.
Suppose $\mu/2$ or more of the singular values of $W_\mu$
are greater than, say, $\varepsilon^{1/3}$.  Then it is
exponentially unlikely that $q^T W_\mu$ has all entries
${\le}\, \varepsilon^{1/2}$, since this would require $\mu/2$
independent random inner products of $q$ with a left singular
vector of $W_\mu$ each to be of order $\varepsilon^{1/6}$,
which is ${\ll}\,1$.  So we may assume that $\mu/2$ of the
singular values of $W_\mu$ are ${\le}\,\varepsilon^{1/3}$.
But it is exponentially unlikely that a random $\mu\times n$
matrix has such a $\mu\times \mu$ submatrix.  This is established
in Lemma~1 of the Appendix.\qed

\section{Conjecture: \boldmath$n^{1/2}$ for the growth factor $\rho$}

I believe that the probabilistic bound of Theorem 2 on the final entry $u_{nn}^{}$
of $U$ is also valid for the growth factor as defined by (\ref{rho}), which
incorporates earlier entries $u_{ij}^{(k)}$ too.
However, a proof has not yet been found, so this is reported as a conjecture.
As mentioned in the Preliminary Note, John Urschel may now have
settled this conjecture in the affirmative~\cite{urschel}.

\medskip

{\em {\bf Conjecture.}  Let GE be applied to a random $n\times
n$ matrix.  It is exponentially unlikely that the growth factor
$\rho$ is much greater than $n^{1/2}$.}

\section{Discussion}

It is remarkable how little interest has been shown in explaining
the rarity of large growth factors in GE, a mystery at the very
heart of numerical analysis.  These quotes from six decades
show how little the subject has moved over the years.

\medskip

{\em Wilkinson\/} 1954 \cite{wil54}:
``It is possible to construct sets [of equations] in which
this factor is achieved but in practice an increase seldom
takes place.''

{\em Wilkinson\/} 1963 \cite{wil63}:
``There are rather special matrices for which the upper bound
is attained$\dots,$ but in practice such growth is very rare.''

{\em Dongarra, Bunch, Moler, and Stewart\/} 1979 \cite{dbms}:
``A great deal of experience has convinced us that the matrices
which lead to the exponentially large coefficients are very
rarely encountered in practice.''

{\em Stoer and Bulirsch\/} 1980 \cite{sb}:
``This bound is much too pessimistic in most cases.''

{\em Gill, Murray and Wright\/} 1991 \cite{gmw}:
``It is rare for the growth factor to exceed 10 on `real'
problems, regardless of the size of $n$.''

{\em Demmel\/} 1997 \cite{dem}:
``The main empirical observation, justified by decades of
experience, is that [$\kern 1.2pt \rho\kern 1.2pt $] is almost
always $n$ or less$\dots$  so we can say that GEPP is {\em
backward stable in practice,} even though we can construct
examples where it fails.''

{\em Higham\/} 2002 \cite{hig}:
``Despite the existence of matrices for which $\rho_n$ is large
with partial pivoting, the growth factor is almost invariably
small in practice.''

{\em Watkins\/} 2010 \cite{wat}:
``Years of testing and experience have shown that the exponential
element growth exhibited by the class of matrices in Exercise
2.7.26 is extremely rare in practice.''

{\em Golub and Van Loan\/} 2013 \cite{gvl}:
``Although there is still more to understand about $\rho$,
the consensus is that serious element growth in Gaussian
elimination with partial pivoting is {\em extremely\/} rare.
{\em The method can be used with confidence.}''

\medskip 

The allusion to `real' problems in the quote by Gill, Murray
and Wright suggests the idea that although instability of
GE may appear for certain matrices, the matrices that arise in
practice somehow avoid this dangerous regime.  I do not believe
such an explanation can work.
A dangerous regime could be avoided by trillions of matrix problems
arising in thousands of applications over a period of half a
century
only by being very small.
In this sense, the stability of GE for random
matrices can be said to {\em explain\/} its stability observed
in practice, though it could never {\em imply\/} stability,
since it might have happened that the unstable matrices arose
exceptionally often in practice.

It would be good to have a mechanistic explanation of the
stability of GE, not just a statistical one.  Here we have nothing
to add\footnote{However, see~\cite{unst}.}
to the discussion in~\cite{ts}, which focusses
on the implications of the
fact that the correction introduced at each step of elimination
is of rank~1.  A key
empirical observation of~\cite{ts} is that for many matrices,
after a few steps of elimination, the distributions of matrix
entries in GE begin to approximate those associated with
$N(0,1)$ initial matrices, suggesting that effects associated
with the central limit theorem are at work, related perhaps to
the ``universality in numerical computations with random data''
investigated by Deift at al.~\cite{dmot}, but we make no claims
here as to the universality of such thermalization in GE or
its relationship to stability in practice.

One interpretation of the present work
might be that the stability of Gaussian
elimination with partial pivoting can be thought of as a law of
numerical computation like the ideal gas law or the second law
of thermodynamics in physics: the justification is statistical,
but overwhelming.  A less positive view might be that we lack
sufficient understanding of the mechanics of GE to make such
an analogy.  Still, it would be odd indeed to argue that the
curves of Figure~\ref{fig-DM} are unrelated to the
success of GE over the years, and now, we are at least part-way to
a theorem that justifies these curves.

\section*{Appendix}

The following lemma, which is related to the mathematics of
restricted isometry constants in the field of compressed
sensing~\cite{ct}, was used in the proof of Theorem~2.
As earlier, ``exponentially unlikely'' is defined by equation
(\ref{highly}).  By a $\mu\times \mu$ submatrix of the $\mu\times
n$ matrix $A$, we mean the matrix formed by selecting some $\mu$
distinct columns of $A$, not necessarily contiguous.

\medskip

{\em {\bf Lemma 1.}
Let\/ $\mu$ be an integer function of\/ $\varepsilon \ll
1 $ with $\mu(\varepsilon) \sim C |\log\varepsilon|$ as
$\varepsilon\to 0$ for some $C>0$ and let $A$ be a random
$\mu\times n$ matrix.  It is exponentially unlikely that $A$
has a\/ $\mu\times \mu$ submatrix with $\mu/2$ or more singular
values ${\le}\,\varepsilon$ with $\varepsilon \ll 1$.}

\medskip

{\em Proof.}
The total number of $\mu\times \mu$ submatrices of $A$ is
${n\choose \mu}$, which is typically about
$(e\kern .4pt n/\mu)^\mu$ and 
can be shown to be always ${\le}\,(3n/\mu)^\mu$.
Rather than needing a finer analysis taking into account
dependencies of the columns of $A$, it is accordingly enough for
us to show that among $(3n/\mu)^\mu$ random
$\mu\times \mu$
matrices, it is exponentially unlikely that there will be
one with $\mu/2$ or more singular values ${\le}\,\varepsilon$ with
$\varepsilon \ll 1$.  (We are employing what is
known as a union bound.)
More precisely, following the definition
(\ref{highly}), it is enough to show that for any $\alpha
>0$ and any $M>0$, with $\varepsilon = n^{-\alpha}$,
for a $\mu\times \mu$ random matrix with $\mu\sim C|\log \varepsilon|$,
the probability that $\mu/2$ or more of its singular values are
${\le}\, \varepsilon$ satisfies
\begin{equation}
p < n^{-M} / (3n/\mu)^\mu \quad \hbox{for all sufficiently large $n$}.
\label{probbound1}
\end{equation}
Since $\mu\to\infty$ as $n\to\infty$, a tighter condition
that it is accordingly sufficient to establish is
\begin{equation}
p < n^{-\mu}  = \varepsilon^{\mu/\alpha}\quad \hbox{for all sufficiently large $n$}.
\label{probbound2}
\end{equation}
Now in fact, by Theorem~1.2 of~\cite{sza}, the probability $p$
satisfies
\begin{equation}
p < C \varepsilon^{\mu^2/4}
\label{probbound3}
\end{equation}
for a universal constant $C$, and this implies 
(\ref{probbound2}).\qed

\end{document}